\input gtmacros
%
\medskipamount=6pt plus 2pt minus 2pt
\mathsurround0.75pt
\hyphenation{handle-body}
\hyphenation{diffeo-morphism}
\hyphenation{diffeo-morphic}
\def \fecsum{Definition\ 2.1} 
\def \fonecomp{Proposition\ 2.5} 
\def \fegAnyA{Example\ 3.6} 
\def \sestsec{4} 
\def \festF{Theorem\ 4.3} 
\def \fsrf{Example\ 4.4} 
\def \fcpt{Example\ 5.10} 
\def \fums{Corollary\ 5.5} 
\def \fgoktso{Example\ 5.6} 
\def \fmsofM{Theorem\ 6.4} 

\catcode`?=11

\def\?empty{}
\long\def\hide#1\endhide{}

\newwrite\labellist
\newread\tmpinput

\newif\ifwritelabellist\global\writelabellistfalse
\font\smallfont=cmr6

\long\def\?loop#1\?repeat{\def\?body{#1}\?iterate}
\long\def\?iterate{\?body\let\?next=\?iterate\fi\?next}


\newbox\arrowtmp
\newbox\arrowtmpA
\newdimen\arrowwd

\def\limitsarrow#1#2#3{%
\setbox\arrowtmp=\hbox{$\scriptstyle11#1$}%
\setbox\arrowtmpA=\hbox{$\scriptstyle11#2$}%
\ifdim\wd\arrowtmp<\wd\arrowtmpA\relax
\arrowwd=\wd\arrowtmpA\relax\else%
\arrowwd=\wd\arrowtmp\relax\fi%
\smash{\mathop{\vbox to 3pt{\vss\hbox to \arrowwd{#3}}}\limits^{#1}\limits_{#2}}%
}

\def\RA#1{
\limitsarrow{#1}{}{\rightarrowfill}}


\def\Rzb#1{\hbox to 0pt{$#1$\hss}}
\def\Lzb#1{\hbox to 0pt{\hss$#1$}}
\def\Czb#1{\hbox to 0pt{\hss$#1$\hss}}
\def\rzb#1{\hbox to 0pt{$\scriptstyle#1$\hss}}
\def\lzb#1{\hbox to 0pt{\hss$\scriptstyle#1$}}
\def\czb#1{\hbox to 0pt{\hss$\scriptstyle#1$\hss}}
\def\ssRzb#1{\hbox to 0pt{$\scriptscriptstyle#1$\hss}}
\def\ssLzb#1{\hbox to 0pt{\hss$\scriptscriptstyle#1$}}
\def\ssCzb#1{\hbox to 0pt{\hss$\scriptscriptstyle#1$\hss}}

\newbox\n?box
\newdimen\n?dimen

\def\Bar#1{\setbox\n?box = \hbox{$#1$}
\kern .08\wd\n?box{\overline{\hbox to .84\wd\n?box{\vphantom {$#1$}\hss}}}
\kern .08\wd\n?box\kern -\wd\n?box\box\n?box}
\def\ctrlBar#1#2#3{\setbox\n?box=\hbox{$#1$}%
\n?dimen=\wd\n?box %
\advance\n?dimen by -#2pt\advance\n?dimen by -#3pt%
\kern #2pt\overline{\hbox to\n?dimen{\hss\vphantom{$#1$}}}%
\kern#3pt\kern-\wd\n?box\box\n?box}

\def\disjointunion{\hbox{$\perp\hskip -4pt\perp$}}



\def\Z{{\Bbb Z}}
\def\Q{{\Bbb Q}}
\def\R{{\Bbb R}}

\newcount\numberscheme
\newcount\label\label=0
\newcount\sectionnumber\sectionnumber=0
\newif\ifexpandthmlabels\expandthmlabelsfalse
\newif\ifexpandbiblabels\expandbiblabelsfalse
\newif\ifexpandsymbollabels\expandsymbollabelsfalse

\numberscheme=1 

\def\check?#1#2{
\xdef\r?r{\expandafter\meaning\csname f?#1\endcsname}
\xdef\rr?r{\expandafter\meaning\csname relax\endcsname}%
\ifx\r?r\rr?r\relax
\xdef\r?r{\expandafter\meaning\csname f#1\endcsname}
{\ifx\r?r\rr?r\relax
\message{#1\space is notdefined.}
XXX%
\else
\edef\x?x{\csname f#1\endcsname ?}\expandafter#2\x?x%
\fi}%
\else
\edef\x?x{ \csname f?#1\endcsname ?}\expandafter#2\x?x%
\fi%
}

\def\F?lakeI#1\ #2?{#1}
\def\F?lakeII#1\ #2?{#2}

\def\partII#1{\check?{#1}{\F?lakeII}}

\def\twothms#1#2{%
\edef\zm?{\csname f?#1\endcsname ?}%
\edef\xx?x{\expandafter\F?lakeI\zm?}%
\edef\zm?{\csname f?#2\endcsname ?}%
\edef\xy?x{\expandafter\F?lakeI\zm?}%
\ifx\xx?x\xy?x\relax\xy?x s \partII{#1} and \partII{#2}\else%
\csname f?#1\endcsname\ and \csname f?#2\endcsname\fi}

\def\stufflabel{\ifnum\numberscheme=1 \number\sectionnumber.\number\label
\else\number\label\fi}

\def\tmplabel{xx}


\def\autonewlabel#1#2#3{\global\advance\label by 1 %
\??autonewlabel{#1}{#2}{#3}}

\def\??autonewlabel#1#2#3{%
\def\xxcc{#1}\ifx\xxcc\tmplabel\relax\else
\xdef\r?r{\expandafter\meaning\csname f?#1\endcsname}
\xdef\rr?r{\expandafter\meaning\csname relax\endcsname}%
\ifx\r?r\rr?r\relax\else%
\xdef\r?r{\csname f?#1\endcsname}
\let\nl?c=\newlinechar%
\newlinechar=`\^^J%
\errmessage{\string f#1\space is already defined as^^J\r?r}%
\let\newlinechar=\nl?c%
\fi\fi%
\def\r?r{#3}%
\ifx\r?r\?empty\relax\edef\r?r{\stufflabel}\fi
\ifexpandthmlabels\relax
\expandafter\xdef\csname f#1\endcsname{\?addlabel{f#1}{#2\ \r?r}}%
\else%
\expandafter\xdef\csname f#1\endcsname{#2\ \r?r}\fi
\expandafter\xdef\csname f?#1\endcsname{#2\ \r?r}
\ifwritelabellist%
\immediate\write\labellist{\def\string\f#1{#2\ \r?r}}
\fi
}

\def\?addlabel#1#2{\setbox0=\hbox{\raise 8pt\hbox to 0pt{%
$\scriptscriptstyle\backslash${\smallfont #1}\hss}{#2}}\unhbox0}

\def\newT#1{\autonewlabel{#1}{\?thm?}{}}
\def\newC#1{\autonewlabel{#1}{\?cor?}{}}
\def\newL#1{\autonewlabel{#1}{\?lem?}{}}
\def\newD#1{\autonewlabel{#1}{\?def?}{}}
\def\newP#1{\autonewlabel{#1}{\?pro?}{}}

\def\newE#1{\autonewlabel{#1}{\?exa?}{}}

\def\newR#1{\autonewlabel{#1}{\?rem?}{}}
\def\newRs#1{\autonewlabel{#1}{\?res?}{}}


\def\?thm?{Theorem}
\def\?cor?{Corollary}
\def\?lem?{Lemma}
\def\?def?{Definition}
\def\?pro?{Proposition}
\def\?for?{Formula}
\def\?exa?{Example}
\def\?dia?{Diagram}
\def\?rem?{Remark}
\def\?res?{Remarks}
\def\ab??X#1\ #2?{\def\x??x{#1}%
\ifx\x??x\?thm? Theorem\else
\ifx\x??x\?cor? Corollary\else
\ifx\x??x\?lem? Lemma\else%
\ifx\x??x\?def? Definition\else%
\ifx\x??x\?pro? Proposition\else%
\ifx\x??x\?for? Formula\else%
\ifx\x??x\?exa? Example\else%
\ifx\x??x\?dia? Diagram\else%
\ifx\x??x\?rem? Remark\else%
\ifx\x??x\?res? Remarks\else%
\x??x%
\fi\fi\fi\fi\fi\fi\fi\fi\fi\fi\ #2}

\def\ab?X#1?{\edef\x??x{#1}\show\x??x\expandafter\ab??X\x??x}
\def\abrev#1{\check?{#1}{\ab??X}}

\newcount\referencecount\referencecount=0
\newif\ifusereferencecount\usereferencecountfalse


\def\ref?list{}

\def\newref#1#2{\global\advance\referencecount by 1 %
\xdef\r?r{\expandafter\meaning\csname b?#1\endcsname}
\xdef\rr?r{\expandafter\meaning\csname relax\endcsname}%
\ifx\r?r\rr?r\relax\else%
\xdef\r?r{\csname b?#1\endcsname}%
\let\nl?c=\newlinechar%
\newlinechar=`\^^J%
\errmessage{\string b#1\space is already defined as^^J\r?r}%
\let\newlinechar=\nl?c%
\fi%
\ifusereferencecount\relax\def\r?r{\number\referencecount}\else%
\def\r?r{#2}\ifx\r?r\?empty\relax\def\r?r{\number\referencecount}%
\fi\fi%
\ifwritelabellist\immediate\write\labellist{%
\def\string\b #1{\r?r}}\fi%
\ifexpandbiblabels\relax%
\expandafter\xdef\csname b#1\endcsname{\?addlabel{b#1}{\r?r}}%
\else%
\expandafter\xdef\csname b#1\endcsname{\r?r}\fi%
\expandafter\xdef\csname b?#1\endcsname{\r?r}%
\expandafter\append\csname b#1\endcsname\to\ref?list%
}

\newtoks\?ta\newtoks\?tb
\long\def\append#1\to#2{\?ta{\\{#1}}\?tb=\expandafter{#2}%
\edef#2{\the\?tb\the\?ta}}
\def\lop#1\to#2{\expandafter\lopoff#1\lopoff#1#2}
\long\def\lopoff\\#1#2\lopoff#3#4{\def#4{#1}\def#3{#2}}

\def\symbolhere#1#2{%
\xdef\r?r{\expandafter\meaning\csname s?#1\endcsname}
\xdef\rr?r{\expandafter\meaning\csname relax\endcsname}%
\ifx\r?r\rr?r\relax\else%
\xdef\r?r{\csname s?#1\endcsname}%
\let\nl?c=\newlinechar%
\newlinechar=`\^^J%
\errmessage{\string #1\space is already defined as^^J\r?r}%
\let\newlinechar=\nl?c%
\fi%
\ifwritelabellist\immediate\write\labellist{%
\def\string\s #1{#2}}\fi%
\ifexpandsymbollabels\relax%
\expandafter\xdef\csname s#1\endcsname{\?addlabel{#1}{#2}}%
\else%
\expandafter\xdef\csname s#1\endcsname{#2}\fi%
\expandafter\xdef\csname s?#1\endcsname{#2}%
}



\def\LproofStyle{%
\def\pf{\ppar\noindent{\bf Proof}\rm\stdspace}
\def\Pf##1{\ppar\noindent{\bf Proof ##1}\rm\stdspace}
}

?rm=cmr8

?Font=cmbx12

\def\smallrm{\small?rm\textfont0=\small?rm\scriptfont0=\scriptscriptfont0%
\textfont1=\scriptfont1\scriptfont1=\scriptscriptfont1%
\textfont2=\scriptfont2\scriptfont2=\scriptscriptfont2%
\let\??xxBBB=\Bbb\def\Bbb{\??xxBBB\scriptstyle}
}

\newbox\addressbox
\newdimen\oodim
\newskip\Titleskip\Titleskip=6pt
\newskip\Aftertopmatterskip\Aftertopmatterskip=10pt
 
\newbox\endbox

\let\??end=\end
\def\end{\box\endbox\vfill\eject\??end}


\def\resetsection{\advance\sectionnumber by 1 %
\ifnum\numberscheme=1 \label=0 \fi}

\def\LsectionStyle{%
\def\newsec##1{\resetsection\vskip 15pt plus 5pt minus 5pt\penalty-1000 %
{\large\bf ##1}\nobreak\vskip 8pt plus 4pt minus 4pt\nobreak}
\def\longnewsec##1##2{\resetsection\bigskip\penalty-1000 %
\currentalignment{\bf ##1}\currentalignment{\bf ##2}
\nobreak\bigskip\nobreak}
}

\def\LciteStyle{%
\def\cite[##1]{{[{##1}]}}
\def\Cite[##1,##2]{[{##1},##2]}
}

\let\?Bbb=\bf

\expandthmlabelstrue
\expandbiblabelstrue
\expandsymbollabelstrue

\def\Lstyle{\LproofStyle\LciteStyle\LsectionStyle}

\def\noLabels{\expandthmlabelsfalse
\expandbiblabelsfalse
\expandsymbollabelsfalse}

\newdimen\labellength
\newdimen\b?w
\newdimen\b??W
\b??W=4pt

\def\bibitem#1{\lop\ref?list\to\?ta
\noindent\llap{[\?ta]\stdspace}}
\catcode`?=12

\Lstyle
\noLabels

\newref{BE}{}
\newref{Cass}{}
\newref{DF}{}
\newref{Ding}{}
\newref{DJDGII}{}
\newref{FAnn}{}
\newref{FJDG}{}
\newref{FTUniv}{}
\newref{Fur}{}
\newref{GompfM}{}
\newref{Gompf}{}
\newref{Gompfcombs}{}
\newref{GrRem}{}
\newref{KirbB}{}
\newref{Kirb}{}
\newref{KS}{}
\newref{Milnor}{}
\newref{Tau}{}
\newref{Wall}{}

\def\ginv#1{b_{{#1}}}

\def\open#1{\setbox0=\hbox{${#1}$}%
\dimen0=\ht0\advance\dimen0 by 2pt%
\raise \dimen0\hbox to\wd0{\hfil\hskip2pt$%
\scriptstyle\circ$\hfil}\hskip-\wd0#1}

\def\ecs{\natural}
\def\mecs#1#2{\natural^{#1}\>{#2}}
\def\ft{{\bf F}_2}
\def\eeight{{\bf E}_8}
\def\hypb{{\bf H}}
\def\calE(#1){{\cal E}(#1)}
\def\insm#1{{\R^4_{#1}}}

\input gtoutput

\volumenumber{1}\volumeyear{1997}\papernumber{6}
\pagenumbers{71}{89}\published{6 December 1997}
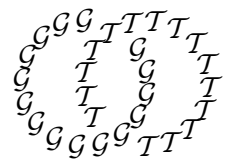
\received{21 October 1997}\accepted{24 November 1997}
\proposed{Robion Kirby}\seconded{Ronald Stern, Ronald Fintushel}

\title{An Invariant of Smooth 4{--}Manifolds}
\shorttitle{An invariant of smooth 4--manifolds}
\author{Laurence R Taylor}
\address{%
Department of Mathematics\\University of Notre Dame\\%
Notre Dame, IN 46556}
\email{taylor.2@nd.edu}

\abstract
We define a diffeomorphism invariant of smooth $4${--}manifolds 
which we can estimate for many smoothings of $\R^4$ 
and other smooth $4${--}manifolds.
Using this invariant we can show that uncountably many smoothings of 
$\R^4$ support no Stein structure.
Gompf [\bGompf] constructed uncountably 
many smoothings of $\R^4$
which do support Stein structures.
Other applications of this invariant are given.
\endabstract
\primaryclass{57R10}\secondaryclass{57N13}\keywords{Smooth 4--manifolds, 
Stein manifolds, covering spaces}

\maketitlepage

\newsec{\number\sectionnumber\quad An invariant of smooth $4${--}manifolds}
\symbolhere{secone}{\number\sectionnumber}

\def\caln{{\cal S}p}
We define an invariant of a smooth $4${--}manifold $M$,   
denoted $\gamma(M)$, by measuring the complexity of
the smooth $\R^4$'s which embed in it.
In our applications, we will have a smooth $4${--}manifold $M$ and
an {\sl exhaustion\/} of $M$ by smooth submanifolds: ie a collection
of smooth submanifolds
$W_0\subset\cdots\subset W_{i-1}\subset\open W_i
\subset W_i
\subset\cdots\subset M$ with $M=\cup W_i$.
Here $\open W_i$ denotes the interior of $W_i$.
We will need to estimate
$\gamma(M)$ from the $\gamma(W_i)$, specifically we want
$\gamma(M)\leq\max_i\bigl\{\gamma(W_i)\bigr\}$.
The definition has several steps.

First, given a smoothing of $\R^4$, $E$, define $\ginv{E}$, 
a non{-}negative integer or $\infty$, as follows.
Let ${\caln}(E)$ denote the collection of smooth, compact, Spin
manifolds $N$ without boundary with hyperbolic intersection form
such that $E$ embeds smoothly in $N$.
Define $\ginv E=\infty$ if ${\caln}(E)=\emptyset$.
Otherwise define $\ginv E$ so that
$$2\ginv E=\min_{N\in {\caln}(E)} \{ \beta_2(N) \}$$
where $\beta_2$ denotes the second Betti number.

While $\ginv E$ measures the complexity of a smooth $\R^4$, 
our ignorance of the properties of smoothings of $\R^4$ prevents
us from using this invariant directly.
For example, we can not rule out the possibility that there is a
smoothing of $\R^4$, $E_\omega$ with the property that all its
compact subsets embed in $S^4$ but $\ginv{E_{\omega}}>0$.
If we were to define $\gamma(M)$ as below using $\ginv{E}$ 
directly, our inequality in the first paragraph would be violated.
Hence we proceed as follows. 
For any smooth $4${--}manifold $M$, let 
$\calE(M)$ denote the set of topological embeddings, 
$e\co D^4\to M$ satisfying two additional conditions. 
First we require $e(\partial D^4)$ to be bicollared and second
we require the existence of a point $p\in\partial D^4$ 
so that $e$ restricted to a neighborhood of $p$ is smooth.
These two conditions are introduced to make the proofs
which follow work more smoothly.
The smooth structure on $M$ induces via $e$ a smooth structure on
the interior of $D^4$, denoted $\insm e$.
Define $\ginv{e}=\ginv{\insm e}$.
Define 
$$\gamma(E)=\max_{e\in{\cal E}(E)}\>\{\>\ginv e\>\}$$ 
If $M$ is Spin, define $\gamma(M)$ to be the maximum of 
$\gamma(E)$ where $E\subset M$ is an open subset homeomorphic 
to $\R^4$.
Note $\gamma(M)$ also is the maximum of $\ginv e$ for $e\in{\cal E}(M)$.
If $M$ is orientable but not Spin and has no compact dual to $w_2$, set 
$\gamma(M)=-\infty$.
If there are compact duals to $w_2$, then define $\gamma(M)$ to be
the maximum of $\gamma(M-F)-\dim_{\ft}H_1(F^2;\ft)$
where $F$ runs over all smooth, compact surfaces in $M$
which are dual to $w_2$.
If $M$ is not orientable, let $\gamma(M)=\gamma(\widetilde M)$
where $\widetilde M$ denotes the orientable double cover of $M$.
The non{-}oriented case will not be mentioned further.

Clearly $\gamma(M)$ is a diffeomorphism invariant.
If $M_2$ is Spin and if $M_1\subset M_2$, or even if every compact 
smooth submanifold of $M_1$ smoothly embeds in $M_2$, then clearly 
$\gamma(M_1)\leq\gamma(M_2)$. 
If $M_2$ is not Spin then such a result is false.
\fcpt\ says $\gamma(CP^2)=0$ and there are
smoothings of Euclidean space, $E$ with $E\subset CP^2$ and
$\gamma(E)$ arbitrarily large.

Note $M$ has no compact duals to $w_2$, if and only if 
$\gamma(M)=-\infty$.
If $M$ has a compact dual to $w_2$, say $F$, then embedding a standard 
smooth  disk in $M-F$ shows that $\gamma(M)\geq-\dim_{\ft}H_1(F;\ft)$.
In particular, if $M$ is Spin, then $\gamma(M)\geq 0$.
The proof of the inequality
$\gamma(M)\leq\max_i\bigl\{\gamma(W_i)\bigr\}$ 
will be left to the reader.
If $M$ is Spin this inequality is an equality.

The precise set of values assumed by $\gamma$ is not known.
One extreme, $-\infty$, is assumed, but no example
with $-\infty<\gamma(M)<0$ is known.
Turning to the non-negative part, $\infty$ is assumed and
it follows from Furuta's work that any non-negative integer
is assumed. 
Adding Taubes's work to Furuta's, 
there are uncountably many distinct smoothings of $\R^4$ 
with $\gamma(E)=n$ for each integer $0\leq n<\infty$.
(The case $n=0$ uses work of  DeMichelis and Freedman [\bDF].)
The referee has noticed that $\infty$ is assumed uncountably often.

We now turn to some applications of the invariant.
In Section \sestsec\ we will estimate $\gamma(M)$ under a condition that is 
implied by the existence of a handlebody structure with no 
$3${--}handles.
It is a theorem that Stein $4${--}manifolds 
have such a handlebody decomposition.
Some version of this theorem goes back to Lefchetz with further
work by Serre and Andreotti--Frankel.
There is an excellent exposition of the Andreotti{--}Frankel theorem
in Milnor [\bMilnor, Section 7 pages 39--40].
It follows from \abrev{estF}\ that if $M$ is a Stein $4${--}manifold,
$\gamma(M)\leq\beta_2(M)$.
In particular, $\gamma(M)=0$ if $M$ is a Stein manifold homeomorphic
to $\R^4$, \abrev{srf}.
Since there are uncountably many smoothings of $\R^4$ with 
$\gamma>0$, \abrev{ums}, there are uncountably many smoothings of
$\R^4$ which support no Stein structure.
In contrast, Gompf \cite[\bGompf] has constructed uncountably many
smoothings of $\R^4$ which do support Stein structures.
They all embed smoothly in the standard $\R^4$ and hence have
$\gamma=0$ as required.
These remarks represent some progress on Problem 4.78 \cite[\bKirb].

The invariant can be used to show some manifolds can not 
be a non-trivial cover of any manifold.
For example, if $E$ is a smoothing of $\R^4$ and is 
a non-trivial cover of some other smooth manifold
then, for $r\geq2$,
$\gamma(\mecs rE)=\gamma(E)$, where $\mecs rE$ denotes
the end-connected sum, \abrev{ecsum}, of $r$ copies of $E$.
By \abrev{goktso}, for each integer $i>0$ there are smoothings $E_i$ 
of $\R^4$ with $\gamma(\mecs rE_i)>2r\gamma(E_i)/3$, so 
these manifolds are not covers of any smooth $4${--}manifold.
Gompf \Cite[\bKirb, Problem 4.79A] asks for smoothings of
$\R^4$ which cover compact smooth manifolds.
These $E_i$ are ruled out, perhaps the first such examples known 
not to cover.

Another use of the invariant is to construct countably many distinct 
smoothings on various non-compact $4${--}manifolds.
The genesis of the idea goes back to Gompf \cite[\bGompfM] 
and has recently been  employed again by Bi\v zaca and Etnyre \cite[\bBE].
In the present incarnation of the idea, one constructs a smoothing, $M$,
with $\gamma(M)=\infty$ and an exhaustion of $M$ by manifolds 
$M(\rho)$, $0\leq \rho<\infty$, homeomorphic to $M$.
One then proves $\gamma(M(\rho_1))\leq\gamma(M(\rho_2))<\infty$
for all $0\leq \rho_1<\rho_2<\infty$ and 
$\lim_{\rho\to\infty}\gamma(M(\rho))=\gamma(M)$.
It follows that $\gamma(M(\rho))$ takes on infinitely many values.
Some hypotheses are necessary: see \abrev{msofM}\ for details.
In particular, \abrev{egAnyA}, there exists a 
Stein $4${--}manifold $M$ with
$H_1(M;\Z)=\Q$ and $H_2(M;\Z)=0$ and this manifold has at
least countably many smoothings.
This may be the first example of a non-compact, 
connected manifold with many smoothings and no topological collar 
structure on any end.
For a non-compact, orientable 4{--}manifold $M$
with no $3$ or $4${--}handles and with finitely generated but infinite
$H^2(M;\Z)$, Gompf has a second construction of countably many
distinct smoothings \Cite[\bGompfcombs, after Theorem 3.1].
The precise relationship between Gompf's smoothings and ours
is not clear.

The author would like the thank R~Kirby, D~Kodokostas and the referee 
for remarks which have substantially improved the exposition and 
scope of the results presented here. The author
was partially supported by the N.S.F.

\newsec{\number\sectionnumber\quad Some properties of the invariant}

The main result in this section is that $\gamma(M)=\gamma(M-K)$
for $K$ a 1{--}complex, \abrev{onecomp} below.
Along the way, we show that $e(D^4)$ can be engulfed in a manifold with 
only zero, one and two handles
for any $e\in{\cal E}(M)$.
These results follow from some nice properties of  
embeddings in $\calE(M)$.

Our first result introduces terminology used later.
If $E$ is a smoothing of $\R^4$ and if $K\subset E$ is compact,
then there exist $e\in\calE(E)$ with $K\subset e(\open D^4)$.
This follows from Freedman's work \cite[\bFJDG] which says that
given any smoothing of $\R^4$, say $E$, there is a homeomorphism
$h\co \R^4\to E$ which is smooth almost everywhere.
In particular, $h$ of any standard ball in $\R^4$ 
is an element of ${\cal E}(E)$.
We say that some property of smoothings of $\R^4$ holds
{\sl for all sufficiently large balls\/} in $E$ provided there is some
compact set $K$ such that for every $e\in\calE(E)$ with 
$K\subset e(\open D^4)$, $\insm e$ is a smoothing of $\R^4$ with
this same property.

To describe our next result, we recall a standard definition.
\newD{ecsum}
\proclaim \fecsum \rm
Given two non-compact smooth $4${--}manifolds $M_1$ and $M_2$, 
define the {\sl end-connected sum\/} of $M_1$ and $M_2$ as follows
\cite[\bGompfM].
First take smooth embeddings of $[0,\infty)$ into each $M_i$.
Thicken  each ray up to a tubular neighborhood.
These tubular neighborhoods are diffeomorphic to $(\R_+^4,\R^3)$,
where $\R^4_+$ denotes the standard half space,
$\bigl\{(x_1, x_2, x_3, x_4)\in\R^4\ \vert\ x_4\geq 0\ \bigr\}$.
The end-connected sum, denoted $M_1\ecs M_2$,
is obtained by removing the interiors of these tubular neighborhoods 
and gluing the two $\R^3$'s together 
by an orientation reversing diffeomorphism.
It is well defined up to diffeomorphism once rays are chosen in $M_1$
and $M_2$, but we suppress this dependence in our notation.
A ray determines an end in $M_i$ and and it follows from
a result of Wall \cite[\bWall] that properly homotopic rays
are isotopic.
Hence the end-connected sum depends only on the proper
homotopy class of the chosen rays.

Elements in $\calE(M)$ behave well with 
respect to end-connected sum.  
We will use the bicollaring of the boundary of elements in
${\cal E}(M)$ several times, but here is the only place that 
the smooth point on the boundary is really useful.
One could have defined the invariant without requiring
the smooth point on the boundary and then use Freedman as above
to modify the proof below.
Of course it also follows from Freedman's work that this other invariant
is the same as the one defined.

\newP{wd}
\proclaim \fwd
Let $M$ be connected and let 
$e_i\in\calE(M)$, $i=1$, \dots, $T$ for some finite integer $T>1$. 
If $e_i(D^4)\cap e_j(D^4)=\emptyset$ for $i\neq j$ then there is an 
element $e\in\calE(M)$ such that $\open D^4_e$ is the 
end-connected sum of all the $\open D^4_{e_i}$.
If there is a surface $F\subset M$ such that $e_i(D^4)\cap F=\emptyset$, 
$i=1$, \dots, $T$,  
then $e$ can be chosen so that $e(D^4)\cap F=\emptyset$.

\pf
One uses the smooth points in $\partial D^4$ to connect the various
$e_i(D^4)$ by thickening up arcs.
See Gompf \cite[\bGompfM] or \Cite[\bKirbB, page 96] for more details.
\qed\ppar

Next we show that elements in ${\cal E}(M)$ can be engulfed
in handlebodies with no $3$ or $4$--handles.
\newP{rfis}
\proclaim \frfis
Let $M$ be a smooth $4${--}manifold and let $e\in\calE(M)$.
Then there exists a smooth compact submanifold,
$V\subset M$ such that $e(D^4)\subset V$ 
and $V$ has a handlebody decomposition
with no $3$ or $4${--}handles.

\pf
Let $W\subset M$ be a smooth, codimension $0$ submanifold with 
$e(D^4)$ in its interior.
We may assume $W$ is connected and hence has a handlebody 
decomposition with only one $0$--handle and no $4$--handles.
Consider the handlebody decomposition beginning with $\partial W$.
Since $\pi_1\bigl(W-e(D^4)\bigr)=\pi_1(W)$, there are 
smoothly embedded compact $1${--}manifolds in $W-e(D^4)$ 
which are homotopic to the cores of the $1${--}handles.
Wall \cite[\bWall] shows how to do an isotopy to bring the cores of the 
$1${--}handles into $W-e(D^4)$.
Let $W_2$ denote the $0$, $1$ and $2$--handles building 
from the other direction.
Then after the isotopy, $e(D^4)$ 
is contained in the interior of $W_2$. \qed

General position arguments now yield

\newP{onecomp}
\proclaim \fonecomp
Let $K\subset M^4$ be a PL proper embedding of a locally-finite 
complex of dimension $\leq 1$.
Then $\gamma(M)=\gamma(M-K)$.

\pf
If $F\subset M$ is a $2${--}manifold, then we can do an isotopy to get 
$K$ and $F$ separated.
Moreover, $F\subset M-K$ is dual to $w_2$ if and only if 
$F$ is dual in $M$ to $w_2$.
Hence we need only consider $F\subset M-K$ dual to $w_2$ in $M$.

Now $\gamma(M)=-\infty$ if and only if 
$\gamma(M-K)=-\infty$.
Hereafter, assume neither is.
Fix $e\in\calE(M-F)$.
Use \frfis\ to get a compact codimension $0$ 
submanifold, $V$, with $e(D^4)\subset \open V$ and $V\subset M-F$ 
such that $V$ has no $3$ or $4${--}handles.
Now do an isotopy to move $V$ off of $K$.
This shows that there is an $e^\prime\in\calE(M-K-F)$ with
$\ginv e=\ginv{e^\prime}$.\qed

\newsec{\number\sectionnumber\quad Few essential $3${--}handles}

We say that a smooth manifold $M$ has {\sl few essential 
$k${--}handles\/} provided $M$ has an exhaustion, 
$W_0\subset W_1
\subset\cdots\ $ where each $W_i$
is a compact, codimension 0, smooth submanifold such that
$H_k(W_{i+1}, W_i; \Z)=0$ for each $i\geq 0$. 

\newRs{hifh}
\proclaim \fhifh \rm
By excision,
$H_k(W_{i+1}, W_i; \Z)=
H_k(\ctrlBar{W_{i+1}-W_i}01, \partial\hskip1pt W_i;\Z)$ so 
if each pair $(\ctrlBar{W_{i+1}-W_i}01, \partial\hskip1pt W_i)$ 
has a handlebody decomposition with no $k${--}handles then
$M$ has few essential $k${--}handles.
We say {\it few} essential $k${--}handles, because $H_k(W_0;\Z)$ 
may be non-zero.

In the next section, we will use a ``few essential $3${--}handles'' 
hypothesis to estimate the invariant $\gamma$.
The key remark needed is the next result.
\newP{SeE}
\proclaim \fSeE
Suppose $M^4$ is smooth, orientable, connected, non-compact and 
has few essential $3${--}handles.
Then 
$\beta_2(W_i)\leq\beta_2(M)$ 
for all $i\geq 0$.

\pf
By the Universal Coefficients Theorem, $H_3(W_i,W_{i-1};\Q)=0$.
(Note by Poincar\'e duality $H_3(W_i,W_{i-1};\Q)=0$ implies
$H_3(W_i,W_{i-1};\Z)=0$ as well, so we get no better result by
only requiring $H_3(W_i,W_{i-1};\Q)=0$ 
in our definition of few essential $3${--}handles.)
It follows from induction and the long exact sequence of
the triple $(W_{i+j},W_{i+j-1},W_i)$ that
$H_3(W_{i+j},W_i;\Q)=0$ for all $i\geq 0$ and $j>0$.
Letting $j$ go to $\infty$, we see $H_3(M,W_i;\Q)=0$.
Hence $H_2(W_i;\Q)\to H_2(M;\Q)$ 
is injective for all $i\geq 0$. \qed\ppar

The rest of this section is devoted to examples.
The first two use \fhifh.

\newE{egProduct}
\proclaim \fegProduct \rm
For any connected, non-compact $3${--}manifold $M^3$,
$M^3\times\R$  has few essential $3${--}handles.
Indeed, find a handlebody decomposition of $M$ with no $3${--}handles
and this gives an evident handlebody decomposition of $M\times\R$
with no $3${--}handles.

\newE{egStein}
\proclaim \fegStein \rm
If $M^4$ is a smooth manifold with a Stein structure, 
then $M$ has few essential $3${--}handles. 
As remarked in Section 1, a theorem of Andreotti and Frankel 
\Cite[\bMilnor, Section 7 pages 39--40] shows 
that there is a proper Morse function, $f\co M\to [0,\infty)$
with no critical points of index $3$ or $4$.

Gompf \cite[\bGompf] proves a partial converse: if there is a proper
Morse function $f\co M\to [0,\infty)$
with no critical points of index $3$ or $4$ then there is a topological
embedding of $M$ inside itself so that the induced smoothing on $M$
supports a Stein structure.
Gompf's result leads to a full characterization of the
homotopy type of Stein manifolds.

\newT{TStein}
\proclaim \fTStein
Any Stein $4${--}manifold has the homotopy type of a countable,
locally-finite CW complex of dimension $\leq 2$.
Conversely, any countable, locally-finite CW complex 
of dimension $\leq 2$ is the homotopy type 
of a Stein $4${--}manifold.

\pf
Any manifold has the homotopy type of the complex built 
from a Morse function, so the Andreotti and Frankel result shows
Stein manifolds have the homotopy type of countable, 
locally-finite CW complexes of dimension $\leq2$.

Conversely, given any countable, locally-finite CW complex 
of dimension $\leq 2$, there is a smooth orientable $4${--}manifold 
with a Morse function with critical points of index $\leq 2$ of the
same homotopy type.
By Gompf, some smoothing of this manifold supports a Stein
structure.\qed\ppar

We can now use standard homotopy-theoretic constructions to
produce examples of manifolds with few essential $3${--}handles.
\newE{egAnyA}
\proclaim \fegAnyA \rm
Let $G$ be any countable group.
Then $G$ has a countable presentation: there are generators $x_g$,
one for each element $g\in G$ and relations $r_{g,h}=x_gx_h(x_{gh})^{-1}$,
one for each element $(g,h)\in G\times G$.
Put the relations in bijection with the positive integers.
Choose new generators $X_{g,i}$, 
one for each $g\in G$ and each integer $i>0$.
Define new relations, $X_{g,i}X_{g,i+1}^{-1}$ and $R_i=
X_{g,i}X_{h,i}(X_{gh,i})^{-1}$ where the $i${--}th relation in the
first presentation is $r_{g,h}$.
Use this second presentation to construct a connected, 
countable, locally-finite  CW complex of dimension $2$ 
with $\pi_1\cong G$.
From \fTStein, there is a Stein manifold with $\pi_1\cong G$.
One can select a set of relations which precisely generate 
the image of the relations in the free abelian quotient of the 
generators.
This gives a new group $\hat G$ and an epimorphism $\hat G\to G$ 
which induces an isomorphism $H_1(\hat G;\Z)\to H_1(G;\Z)$.
The resulting $2${--}complex has $H_2=0$.
Hence there is a Stein manifold $M$ with $H_1(M;\Z)$ any 
countable abelian group and $H_2(M;\Z)=0$.

Here are two ways to construct additional examples. 
\newE{egDiffatinfinity}
\proclaim \fegDiffatinfinity \rm
If $M_1$ is diffeomorphic at infinity to $M_2$ 
and if $M_2$ has few essential $3${--}handles, 
then $M_1$ has few essential $3${--}handles.
If $M_1$ and $M_2$ have few essential $3${--}handles, 
then  $M_1\ecs M_2$ has few essential $3${--}handles.

\newsec{\number\sectionnumber\quad Some estimates of $\gamma(M)$}
\symbolhere{estsec}{\number\sectionnumber}

\newD{xx}
\proclaim \fxx \rm
An oriented $4${--}manifold has an intersection form on $H_2(M;\Z)$.
We say $M$ is {\sl odd\/} if there is an element whose intersection with
itself is odd: otherwise we say $M$ is {\sl even}.
Spin implies even, but there are manifolds like the Enriques surface
which are not Spin but are even.

\newT{estS}
\proclaim \festS
Let $M$ be the interior of an orientable, smooth, 
compact manifold with boundary (which may be empty).
Then
$$\gamma(M)\leq \beta_2(M)-\cases{0&$M$ even\cr1&$M$ odd\thinspace.\cr}$$

\pf
It suffices to deal with each component of $M$ separately 
so assume $M$ is connected. 
If $\partial M=\emptyset$ replace $M$ by $M-D^4$ where $D^4$
denotes a smooth standard disk.
By \fonecomp, $\gamma(M)=\gamma(M-D^4)$ 
so hereafter assume $\partial M\neq\emptyset$.
Since $M$ is compact, there are compact duals to $w_2$.
Let $F$ be a fixed dual to $w_2$, and then fix an $e\in\calE(M - F)$.
Note it suffices to prove
$\ginv e-d_1(F)\leq\beta_2(M)$, $M$ even, 
(or $\beta_2(M)-1$, $M$ odd) where 
$d_1(F)=\dim_{\ft}H_1(F;\ft)$.
Let $\Bar M$ denote the compact manifold with boundary 
whose interior is $M$.
One can add $1${--}handles to $F$ inside $M-e(D^4)$ to ensure $F$
is connected without changing $d_1(F)$ so assume $F$ is connected.
Let $U=\Bar M - F$ and let $N$ be the double of $U$. 
Since $U$ is Spin, $N$
is a closed, compact Spin manifold with signature 0.
The composition $r\co N\subset U\times[0,1]\ \RA{proj.}\  U$
splits the inclusion, so
$$H_\ast(N;\Q)=H_\ast(U;\Q)\oplus H_\ast(U, \partial U;\Q)$$
since $(U,\partial U)\to (N,U)$ is an excision map.
Hence $\ginv e\leq\ell$ where 
$$\eqalign{%
2\ell+2=&\chi(N)+2\beta_1(N)=
2\cdot\bigl(\chi(M)-\chi(F)+\beta_1(N)\bigr)\cr=&
2\cdot\bigl(\chi(M)+d_1(F)
+\beta_1(N)-2\bigr)}$$
so $\ginv e-d_1(F)\leq \chi(M)+\beta_1(N)-3$.

We now proceed to compute $$H_1(N;\Q)=
H_1(U;\Q)\oplus H_1(U,\partial U;\Q)=H_1(U;\Q)\oplus H^3(U;\Q).$$
Denote the disk{--}sphere bundle pair
for the normal bundle of $F\subset M$ by $\bigl(D(F),S(F)\bigr)$.
The next two sequences are exact:
$$\leqalignno{%
&H^3(\Bar M)\to H^3(U)\to H^4(\Bar M,U)\cong H^4\bigl(D(F),S(F)\bigr)
\cong H_0(F)\cr
&H^2(F)\cong H_2\bigl(D(F),S(F)\bigr)\cong H_2(\Bar M,U)\to
H_1(U)\to H_1(\Bar M)&(\ast\ast)\cr
}$$
Since $H_0(F)=\Q$ and $H^2(F)=\Q$ or $0$,
$\beta_1(U)\leq\beta_1(M)+1$ and $\beta_3(U)\leq\beta_3(M)+1$.
It follows that
$\ginv e-d_1(F)\leq \chi(M)+\beta_1(M)+\beta_3(M)-1=\beta_2(M)$.

If the surface dual to $w_2$ is not orientable, then, in $(\ast\ast)$ above, 
$H^2(F)=0$, so $\beta_1(U)\leq\beta_1(M)$ and 
$\ginv e-d_1(F)\leq \beta_2(M)-1$.
If $F$ is orientable, then $H_2(F)=H_2(\Bar M,U)=\Q$. 
Since the form is odd, there is an element $x\in H_2(\Bar M)$ of odd
self-intersection.
It follows that the image of $x$ in $H_2(\Bar M,U)$ 
is non-zero, so again 
$\beta_1(U)\leq\beta_1(M)$ and 
$\ginv e-d_1(F)\leq \beta_2(M)-1$.\qed

\newT{estF}
\proclaim \festF
If $M$ is orientable with few essential $3${--}handles, then
$$\gamma(M)\leq \beta_2(M)-
\cases{0&$M$ even\cr1&$M$ odd\thinspace.\cr}$$

\pf
If there are no compact duals to $w_2$ the required inequality
is clear, so assume we have compact duals.
Fix a compact dual $F$ and an $e\in\calE(M-F)$ 
and show $\ginv e-d_1(F)\leq\beta_2(M)$, $M$ even
($\beta_2(M)-1$ $M$ odd).
If $\beta_2(M)=\infty$, the required inequality is  
immediate, so assume it is finite.
Use \fSeE\ to find a compact smooth manifold with boundary, 
say $W^4\subset M$,
with $F\disjointunion e(D^4)\subset \open W$ and
$\beta_2(W)\leq \beta_2(M)$.
Then $F$ is dual to $w_2(W)$ so $\ginv e-d_1(F)\leq\gamma(W)$ and
\festS\ implies $\gamma(W)\leq\beta_2(W)$ and 
since $\beta_2(W)\leq \beta_2(M)$, $\gamma(W)\leq\beta_2(M)$.
If the form on $M$ is odd, then we can choose $W$ so large
that the form on $W$ is also odd and then
$\gamma(W)\leq \beta_2(M)-1$.\qed

\newE{srf}
\proclaim \fsrf \rm
If $M$ is a smoothing of $\R^4$ with few essential 
$3${--}handles, then $\gamma(M)=0$.
Hence, $\gamma(M)=0$ if $M$ supports a Stein structure or 
if $M$ is diffeomorphic to $N^3\times\R$.
In particular, the standard smoothing of $\R^4$ has $\gamma=0$ 
and hence any smoothing of $\R^4$ which embeds in it also has 
$\gamma=0$.

\newR{nTH}
\proclaim \fnTH \rm
In the next section, we show many smoothings of $\R^4$ have $\gamma>0$.
It follows that these smoothings must have infinitely many 
$3${--}handles in {\sl any\/} handlebody decomposition.
Even more must be true.
For any exhaustion of such a smoothing of $\R^4$,
$H_3(W_{i+1},W_i;\Z)$ must be non-zero for infinitely many $i$.

\newsec{\number\sectionnumber\quad Examples of $\gamma$
for smoothings of $\R^4$}
\symbolhere{rfsec}{\number\sectionnumber}

First observe that for a smoothing of $\R^4$, $\gamma$ 
depends only on the behavior at infinity.

\newT{ifoi}
\proclaim \fifoi
Let $E_1$ and $E_2$ be smoothings of $\R^4$ and suppose
$E_1$ embeds at $\infty$ in $E_2$ (ie there is an open subset of
$E_2$ which is diffeomorphic at $\infty$ to $E_1$).
Then $\gamma(E_1)\leq\gamma(E_2)$.

\pf
If $E_3\subset E_2$ is the submanifold which is 
diffeomorphic at $\infty$ to $E_1$,
it follows that $E_3$ is homeomorphic to $\R^4$.
It suffices to prove $\gamma(E_1)=\gamma(E_3)$.

Let $F\co E_1-V_1\to E_3-V_2$ be a representative of the 
diffeomorphism at $\infty$.
Let $e_i\in\calE(E_1)$ be a sequence so that 
$e_i(D^4)$ form an exhaustion and further assume
$V_1\subset e_0(\open D^4)$.
Let $W_i\subset E_3$ be $F\bigl( e_i(D^4)-V_1\bigr)\cup V_2$.
The $W_i$ are an exhaustion of $E_3$ and $W_i=\hat e_i(D^4)$ 
for elements $\hat e_i\in\calE(E_3)$.
Check $(\insm{e_i}-V_1)\cup V_2=\insm{\hat e_i}$.
Next check $\ginv{e_i}=\ginv{\hat e_i}$:
if $\insm{e_i}$ embeds in $N$, then
$\insm{\hat e_i}$ embeds in $N^\prime=(N-V_1)\cup V_2$ and
$N^\prime$ is homeomorphic to $N$.\qed\ppar

It is easy to estimate $\gamma$ for an end-connected sum.
\newL{ecsL}
\proclaim \fecsL For $E_1$ and $E_2$ any two smoothings of $\R^4$,
$$\max\bigl\{\gamma(E_1),\gamma(E_2)\bigr\}\leq
 \gamma(E_1\ecs E_2)\leq \gamma(E_1)+\gamma(E_2)\ .$$

\pf
The lower bound follows since $E_i\subset E_1\ecs E_2$.
To see the upper bound, consider the ordinary connected sum,
$E_1\# E_2$.
One can embed a smooth $\R^1$ meeting the $S^3$ in the
connected sum transversely in one point so that 
$E_1\ecs E_2$ and $E_1\# E_2-\R^1$ are diffeomorphic.
By \fonecomp, $\gamma(E_1\ecs E_2)=\gamma(E_1\# E_2)$.
Now, given any embedding $e\in\calE(E_1\# E_2)$, 
there are embeddings 
$e_i\in\calE(E_i)$, each of which contains the disk used to form 
the connected sum and so that $e_1(D^4)\# e_2(D^4)$ contains 
$e(D^4)$ in its interior.
If each $e_i(D^4)$ embeds in $N_i$, $e_1(D^4)\# e_2(D^4)$ embeds 
in $N_1\# N_2$.\qed\ppar

Donaldson \cite[\bDJDGII] began work on the $11/8${--}th's conjecture
(Problem 4.92 \cite[\bKirb]), one version of which says that any 
smooth Spin $4${--}manifold, $N$, with form $2s\eeight\perp t\hypb$ 
must have $t\geq 3s$.
Here $\eeight$ is the even, unimodular form of rank $8$ and index $-8$ 
and $\hypb$ is the rank $2$ hyperbolic form.
While this conjecture is still unsolved in general
Donaldson \cite[\bDJDGII] proved $t>0$ and Furuta \cite[\bFur]
has shown that $t>2s$.
Recall that both these results require no condition on $\pi_1$,
which accounts for our lack of $\pi_1$ conditions on the manifolds
in $\caln(E)$ from Section 1.

We say that a smoothing $E$ of $\R^4$ is {\sl semi-definite\/}
provided there is a positive integer $k$ and a compact, closed, 
topological Spin $4${--}manifold $M$ with form $2s\eeight$, $s>0$, 
so that some smoothing of
$M-pt$ is diffeomorphic at $\infty$ to $\mecs kE$.
If the $k$ and $s$ are important we will say $E$ is 
$(k,s)${--}semi-definite.
If $E_1$ embeds at $\infty$ in $E_2$ and if $E_1$ is 
$(k,s)${--}semi-definite, then $E_2$ is $(k,s)${--}semi-definite.
If $E$ is $(k,s)${--}semi-definite, then all 
sufficiently large balls in $E$ are $(k,s)${--}semi-definite.
We say $E$ is $(k,s)${--}simple-semi-definite if the $M$ can
be chosen to be simply-connected.

\newT{sdrf}
\proclaim \fsdrf
If $E$ is $(k,s)${--}semi-definite, then $\gamma(E)>2s/k$ and 
$$\displaystyle\lim_{r\to\infty}\gamma(\mecs rE)=\infty.$$

\pf
If $\gamma(E)=\infty$ the result is clear, so assume $\gamma(E)<\infty$.
Select a large topological ball in $\mecs kE$ so that the smoothing on
its interior, $E_1$, is also diffeomorphic at $\infty$ to a smoothing 
of $M-pt$, say $V$.
By selecting the ball large enough, $\gamma(\mecs kE)=\gamma(E_1)$.
By definition, $E_1$ embeds  in a smooth, closed, compact, Spin manifold, 
$N^4$, whose intersection form is $\gamma(E_1)\hypb$.
Choose a compact set $K\subset V$ and a ball $\Delta\subset E_1$
so that $V-K$ and $E_1-\Delta$ are diffeomorphic and use the 
diffeomorphism to glue $N-\Delta$ and $V$ together along
$V-K$.
The resulting manifold is smooth and has form 
$2s\eeight\perp \gamma(E_1) \hypb$.
The case $\gamma(E_1)=0$ is forbidden by Donaldson \cite[\bDJDGII]
and by \fecsL, $k\gamma(E)\geq\gamma(E_1)>0$.

By Furuta \cite[\bFur] it further follows that $\gamma(\mecs kE)>2s$.
Since $\mecs{\ell k}E$ is diffeomorphic at $\infty$ to a smoothing of
$M_\ell-pt$ where $M_\ell$ is the connected sum of $\ell$ copies of $M$,
Furuta's result implies $\gamma(\mecs {k\ell}E)>2s\ell$.
Since $\gamma(\mecs rE)$ is a non-decreasing function of $r$,
the limit exists and is $\infty$. \qed\ppar

The next result follows from work of Taubes \cite[\bTau].

\newT{xx}
\proclaim \fxx
If $E$ is a simple-semi-definite
smoothing of $\R^4$ then any 
sufficiently large ball in $E$ is not diffeomorphic to any
larger ball containing it in its interior.

\pf
If not, then for any compact set $K\subset E$ there exists a 
pair of balls $e_i\in\calE(E)$ with 
$K\subset e_1(\open D^4)\subset e_1( D^4)\subset e_2(\open D^4)$
and with $\insm{e_1}$ diffeomorphic to $\insm{e_2}$.
We may further require that $\insm{e_1}$
is $(k,s)${--}semi-definite, 
where $E$ is $(k,s)${--}semi-definite.
We can take the end-connected sum in such a way that
$\mecs k{\insm{e_1}}\subset \mecs k{\insm{e_2}}$
Hence there is a smoothing of $M-pt$ which is diffeomorphic
at $\infty$ to $\mecs k{\insm{e_1}}$, where the form on $M$
is $2s\eeight$. 
But this is forbidden by \cite[\bTau].
\qed

\newC{ums}
\proclaim \fums
If $E$ is a simple-semi-definite smoothing of $\R^4$ with 
$\gamma(E)<\infty$
then there are uncountably many 
simple-semi-definite smoothings of $\R^4$ 
with $\gamma$ equal to $\gamma(E)$.

\pf
All sufficiently large balls inside of $E$ are simple-semi-definite
and have $\gamma$ equal to $\gamma(E)$.\qed

\newE{goktso}
\proclaim \fgoktso
There exists an $e\in\calE(S^2\times S^2)$ so that
$\insm e$ is a $(3,1)${--}simple-semi-definite smoothing of $\R^4$ 
with $\gamma=1$.
Hence $2r/3<\gamma(\mecs r{\insm e})\leq r$.
For each integer $n\geq 1$ there exists an $r_n$ such that
$\gamma(\mecs {r_n}{\insm e})=n$.

\pf
Consider the standard $D^4$ with a Hopf link in its boundary.
Attach two Casson handles with $0$ framing to this Hopf link 
and call the interior a {\sl Casson wedge}.
Let $M^4$ denote the Kummer surface or any other simply-connected
smooth Spin manifold with form $2\eeight\perp 3\hypb$.

Casson's results \cite[\bCass] allow us to construct a particular
Casson wedge, $C$, so that 
$C\disjointunion C\disjointunion C\subset M$.
We can also find $C\subset S^2\times S^2$.
Freedman \cite[\bFJDG] constructs a topological embedding of
$S^2\vee S^2\subset C$ and shows that $E=S^2\times S^2-S^2\vee S^2$
is homeomorphic to $\R^4$.

By construction $\gamma(E)\leq 1$ and 
$\mecs 3E$ is diffeomorphic at $\infty$ to a smoothing of $M^\prime-pt$,
where $M^\prime$ is the simply-connected topological manifold
with intersection form $E_8\perp E_8$ constructed by Freedman.
By \fsdrf, $\gamma(E)>0$.
Now choose $e\in\calE(S^2\times S^2)$ so that $e(D^4)\subset E$
is so large that $\gamma(\insm e)=\gamma(E)=1$ and $\insm e$ 
is $(3,1)${--}semi-definite.

Since the set
$\bigl\{ \gamma(\mecs r{\insm e})\bigr\}$ is unbounded
and  since $0\leq\gamma(\mecs {r+1}{\insm e})-
\gamma(\mecs r{\insm e}))\leq 1$ by \fecsL, 
every integer $n\geq 1$ is
assumed by some $\gamma(\mecs r{\insm e})$.
\qed

\newR{xx}
\proclaim \fxx \rm
If $E$ is any $(3,1)${--}semi-definite smoothing of $\R^4$ with 
$\gamma(E)=1$, it follows as above that for $n\geq1$ there is an
$r_n$ such that $n=\gamma(\mecs{r_n}{E})$.
Moreover, $r_n=n$ for $1\leq n\leq3$.
If the $11/8$-th's Conjecture holds, then 
$r_n=n$ for all $n\geq1$.

Let $U$ denote the interior of the universal half-space constructed
in \cite[\bFTUniv] or any other smoothing of $\R^4$ into which all 
others embed.
Then clearly we have:
\newC{urf}
\proclaim \furf
$\gamma(U)=\infty$.\rm

The above results give some information for compact manifolds.
\newE{xx}
\proclaim \fxx \rm 
If $M$ is the connected sum of $s$ copies of $S^2\times S^2$, then
$2s/3<\gamma(M)\leq s$ and $\gamma(M)=s$ if $s\leq 3$
({\fgoktso}).

\newE{cpt}
\proclaim \fcpt \rm
The case $M=CP^2$ is more interesting.
It is clear $\gamma(CP^2)=0$: the upper bound comes from 
\festS\ while the lower bound comes from a smooth disk 
missing the $CP^1$ dual to $w_2$.
More interestingly, let $M_s$ be the connected sum of $2s$
copies of Freedman's $\eeight$ manifold \cite[\bFJDG]. 
Then there is a homeomorphism
$h\co CP^2\# M_s\to CP^2\#16s\Bar{CP}^2$.
Apply $h$ to the $S^2=CP^1\subset CP^2$ to see a 
locally flat topological embedding 
$\iota\co S^2\to CP^2\#16s\Bar{CP}^2$ 
whose complement is $M_s$.
A neighborhood of $\iota(S^2)$ is homeomorphic to $CP^2-pt$.
With a bit more care, we can choose $\iota$ and a topological
embedding $\iota_0\co S^2\to CP^2$ so that
$\iota$ and $\iota_0$ have diffeomorphic neighborhoods.
Hence $E_s=CP^2-\iota_0(S^2)$ is homeomorphic to $\R^4$
and diffeomorphic at $\infty$ to a smoothing of $M_s-pt$.
See Kirby \Cite[\bKirbB, page 95] for more details.
By construction, $E_s$ is $(1,s)${--}simple-semi-definite, so
\fsdrf\ shows $\gamma(E_s)\geq 2s$.
Hence there are $e\in\calE(CP^2)$ with $\ginv e$ arbitrarily large.
The above calculation of $\gamma$ shows that a dual to $w_2$ in
the complement of the image of $e$ must have large $H_1$.

As pointed out by the referee, more mileage is available from
this example.
With a bit of care, it can be arranged so that 
$E_s\subset E_{s+1}\subset CP^2$.
If we let $E_\infty\subset CP^2$ denote the union, 
then $\gamma(E_\infty)=\infty$ and again with care,
$E_\infty$ can be extended to have a smooth patch of boundary $\R^3$.
Fix an orientation reversing homeomorphism between $\R^4$ and 
some simple-semi-definite smoothing of $\R^4$ and let
$E(\rho)$ denote the smoothing inherited by the ball of radius $\rho$.
Following Ding's use of Taubes's periodic ends theorem,
\cite[\bDing], note that for sufficiently large $\rho$, 
$E_\infty\ecs E(\rho_1)$ is diffeomorphic to 
$E_\infty\ecs E(\rho_2)$ if and only if $\rho_1=\rho_2$.
Since $\gamma(E_\infty\ecs E)=\infty$, there are uncountably
many distinct smoothings of $\R^4$ with $\gamma=\infty$.
Together with \partII{ums}\ and 
\partII{goktso}\ this shows that for each 
$0\leq n\leq\infty$ there are uncountably many
distinct smoothings of $\R^4$ with $\gamma=n$.

\newsec{\number\sectionnumber\quad Smoothings of other $4${--}manifolds}

In this section we will start with a non-compact smooth manifold $M$, 
construct a family of smoothings $M(\rho)$ 
and prove that $\gamma\bigl(M(\rho)\bigr)$
takes on countably many distinct values.
We will start with a smoothing of $M$ satisfying some hypotheses 
which will be clarified later.
Then we form $M\ecs U$, where $U$ is a smooth $\R^4$ into which
all others embed.
Without much trouble, we see $\gamma(M\ecs U)=\infty$.
For positive real numbers $\rho$ we construct a family of 
open submanifolds, $M(\rho)\subset M\ecs U=M(\infty)$
which exhaust $M(\infty)$.
Again without much trouble, we show that
$\gamma\bigl(M(\rho_1)\bigr)\leq \gamma\bigl(M(\rho_2)\bigr)$
for $\rho_1<\rho_2$ and that 
$\gamma\bigl(M(\infty)\bigr)=
\displaystyle\max_{\rho}\bigl\{\>\gamma\bigl(M(\rho)\bigr)\>\bigr\}$.
Finally (the hardest part) we show $\gamma\bigl(M(\rho)\bigr)<\infty$
for $\rho<\infty$.
At this point, it is easy to show that there is a countable collection of
$\rho_1<\cdots < \rho_i<\cdots <\infty$ so that
$\gamma\bigl(M(\rho_{1}\bigr)<\cdots <\gamma\bigl(M(\rho_i)\bigr)
<\cdots<\infty$ so 
we have our countable collection of distinct smoothings.

To show $\gamma\bigl(M(\rho)\bigr)<\infty$ will require constructing
a $4${--}manifold with few essential $3${--}handles.
To be able to do this, we need some structure to start with on $M$.
We actually need our structure on $M\ecs U$ and this manifold never
has few essential $3${--}handles so we need some structure close to
few essential $3${--}handles but weaker.
Roughly, we require an exhaustion of $M$ by topological manifolds
$W_i$ with $H_3(W_{i+1},W_i;\Z)=0$ so that the embedding
$\partial W_i\to M$ is smooth except for a finite number of $3${--}disks.
The precise condition follows next.

We say that an exhaustion of $M$ by topological manifolds $W_i$ is
{\sl almost-smooth\/} provided there exists 
a space $\Delta=\disjointunion_{k=0}^T D^3$ and a proper
embedding
$\tau\co \Delta\times[-1,\infty)\to M$
which is smooth near $(\partial \Delta)\times[-1,\infty)$
and a function $\mu\co M\to[-4,\infty)$ satisfying the following.
Let 
$Y=M-\tau\bigl(\open\Delta\times(-1,\infty)\bigr)-
\tau(\Delta\times-1)$
and note $Y$ is a smooth manifold with boundary 
$(\partial\Delta)\times(-1,\infty)$.
We require
\def\ASone{1}
\def\AStwo{2}
\def\ASthree{3}
\items\item{(\ASone)}%
$W_i=\mu^{-1}\bigl([-4,i]\bigr)$, 
\item{(\AStwo)}$\mu\circ\tau\co\Delta\times[-1,\infty)\to[-4,\infty)$
is projection followed by inclusion,
\item{(\ASthree)}$\mu$ restricted to $Y$ is smooth and the integers in
$(-1,\infty)$ are regular values.\enditems

Note that $\open Y=M-\tau\bigl(\Delta\times[-1,\infty)\bigr)$
is homeomorphic to $M$.
Indeed $M(\rho)$ will be $M-\tau\bigl(\Delta\times[\rho,\infty)\bigr)$
for some almost-smooth exhaustion.
The replacement for ``few essential $3${--}handles'' will be an
almost-smooth exhaustion with $H_3(W_{i+1},W_i;\Z)=0$.

Our first requirement on our structure follows.
\newR{asecs}
\proclaim \fasecs \rm
If $M$ has an almost-smooth exhaustion with 
$H_3(W_{i+1},W_i;\Z)=0$, then
so does $M\ecs E$ for any $E$ homeomorphic to $\R^4$.
Just add one more $3${--}disk to $\Delta$ and 
construct $\tau$ and $\mu$ so that the end-connect sum
takes place inside the added cylinder.

Before introducing the main theorem, we discuss some examples.
\newE{xx}
\proclaim \fxx \rm
If $(N,\partial N)$ is a topological manifold, then $\open N$ has a 
smoothing with an almost-smooth exhaustion with 
$H_\ast(W_{i+1},W_i;\Z)=0$ for all $\ast$.
Indeed, one can smooth $\open N$ so that there is a topological
collar which is smooth off a cylinder $D^3\times[-1,\infty)$.
Here are the details.
Smooth $N-pt$ and choose a smooth collar 
$c\co \partial N\times[-2,\infty]\to N$
so that $c(\partial N\times\infty)=\partial N$.
Let $\mu\co N\to[-4,\infty]$ be
chosen to be smooth on $\mu^{-1}\bigl((-3,\infty]\bigr)$ 
and projection on the collar.
Pick a smooth $\iota\co \Delta=D^3\subset \partial N$ and let
$\tau\co \Delta\times[-1,\infty)
\RA{\iota\times id}\,
\partial N\times[-1,\infty)\RA{c}\open N$.
Note that $\mu^{-1}\bigl([-4,-2])$ is a compact submanifold of $\open N$
with a smooth boundary in the smooth structure on $\open N-pt$;
note further that $pt$ is in its interior.
Let $X=\mu^{-1}\bigl([-4,-2])\cup 
c\bigl(\iota(D^3)\times[-2,\infty)\bigr)$ and note 
$\open N-\open X\subset \open N-pt$ is a smooth submanifold.
Let $M$ be the smoothing on $\open N$ given by smoothing $X$ 
rel boundary and extending over $\open N-\open X$.
As things now stand, $\mu$ may not be smooth on the space
$Y$ in the definition, but it can clearly be altered to be smooth
while leaving it fixed on $\mu^{-1}\bigl[-1,\infty)\bigr)$.
Clearly (\AStwo) and (\ASthree) are satisfied and
$\ctrlBar{W_{i+1}-W_i}02=\partial N\times[i,i+1]$ so 
$H_\ast(W_{i+1},W_i;\Z)=0$.

Here is a way to produce new examples.
\newE{ksi}
\proclaim \fksi \rm
By Kirby and Siebenmann \cite[\bKS] the stable isotopy classes\break 
of smoothings of $M$ are in one{-}to{-}one 
correspondence with elements $\alpha\in$\break$H^3(M;\ft)$.
To construct such a smoothing, select a proper embedding 
$\R^1\to M$ which is dual to $\alpha$.
Extend to an embedding $\hat\alpha\co D^3\times\R^1\to M$.
Freedman \cite[\bFAnn] has constructed a smooth proper homotopy
$D^3\times\R^1$, $(D^3\times\R^1)_\Sigma$, and \cite[\bFJDG] a 
homeomorphism $h\co (D^3\times\R^1)_\Sigma\to D^3\times\R^1$
which is the identity on the boundary and represents the non-zero 
element in $H^3(D^3\times\R^1;\ft)$.
Remove the image of $\hat\alpha$ and replace it by 
$(D^3\times\R^1)_\Sigma$ to get $M_{\hat\alpha}$ and let
$h_\alpha\co M_{\hat\alpha}\to M$ be the evident homeomorphism.
Now let $M$ have an almost-smooth exhaustion, $\tau$ and $\mu$.
The two ends of $\R^1$ and $\hat\alpha$ determine two ends of $M$ 
(which may be the same).
We can further guarantee that the image of $\hat\alpha$ misses the image
of $\tau$.
Let $\Delta^\prime$ be $\Delta$ with two more $3${--}disks added.
One can choose 
$\tau^\prime\co \Delta^\prime\times[-1,\infty)\to M_{\hat\alpha}$ 
and $\mu^\prime$ so that $\tau^\prime$ and $\mu^\prime$ give 
an almost-smooth exhaustion for of $M_{\hat\alpha}$ with the same
topological manifolds $W_i$.

\newT{msofM}
\proclaim \fmsofM
Let $M$ be a non-compact, orientable $4${--}manifold with an
almost-smooth exhaustion with $H_3(W_{i+1},W_i;\Z)=0$.\nl
If $M$ satisfies $\dim_{\ft}H_2(M;\ft)+\beta_2(M)<\infty$, 
then there are at least countably many
distinct smoothings of $M$ in each stable isotopy class.

\pf
By \partII{ksi}\ there is a smoothing of $M$ in a given stable isotopy
class which possesses an almost-smooth exhaustion still satisfying
$H_3(W_{i+1},W_i;\Z)=0$.
Use $M$ to denote this smoothing.
If $U$ is any $\R^4$ into which all others smoothly embed,
form $M\ecs U$. 
This is a smoothing of $M$ in the same stable isotopy class and
has an almost-smooth exhaustion $H_3(W_{i+1},W_i;\Z)=0$ 
by \partII{asecs}.
Fix $\tau$ and $\mu$ to give an almost-smooth
exhaustion of $M\ecs U$ with $H_3(W_{i+1},W_i;\Z)=0$.

For each $\rho$, $0\leq\rho<\infty$, let 
$$M(\rho)=M\ecs U -\tau\bigl(\Delta\times [\rho,\infty)\bigr)$$
and observe each $M(\rho)$ is homeomorphic to $M$ and since
$M(\rho)$ is an open subset of $M\ecs U$ it is a smooth manifold.
Let $M(\infty)=M\ecs U$.
Note that if $\rho_1<\rho_2$, $M(\rho_1)\subset M(\rho_2)$.

Since $\dim_{\ft}H_2(M;\ft)<\infty$, all the $M(\rho)$ have compact duals
to $w_2$ and if $F\subset M(\rho_1)$ is a compact dual to $w_2$,
then it is a compact dual to $w_2$ in $M(\rho_2)$ 
whenever $\rho_1<\rho_2$.
Hence $\gamma\bigl(M(\rho)\bigr)$ is a non-decreasing function 
of $\rho$.
Fix a compact dual to $w_2$ in $M(0)$, say $F_0$.
Given any integer $S$, there exist balls $e\in\calE(U)$ with $\ginv e>S$ 
and there exists a $T$ such that the image of $e$ lies in $M(\rho)$ for
all $\rho>T$.
For $\rho>T$, $\gamma\bigl(M(\rho)\bigr)>S-\dim_{\ft}H_1(F_0;\ft)$
so $\lim_{\rho\to\infty}\gamma\bigl(M(\rho)\bigr)=\infty$.

It remains to prove $\gamma\bigl(M(\rho)\bigr)<\infty$ if $\rho<\infty$.
To do this, we construct a manifold $X(\rho)$ so that $X(\rho)$ has
few essential $3${--}handles.
We show that $\gamma\bigl(M(\rho)\bigr)\leq\gamma\bigl(X(\rho)\bigr)$ 
and that $\beta_2\bigl(X(\rho)\bigr)<\infty$.
\festF\ implies that $\gamma\bigl(X(\rho)\bigr)<\infty$.

To construct $X(\rho)$ first note that we can deform $\mu$ to
a function $\mu_X\co M\ecs U\to [-4,\infty)$
which is smooth everywhere, equal to $\mu$ on $Y$ and has
the non-negative integers as regular values.

Let $k_\rho$ denote the first integer so that 
$\mu^{-1}\bigl([-4,\rho]\bigr)\subset \mu_X^{-1}\bigl([-4,k_\rho]\bigr)$.
Let $Z(\rho)=\tau(\Delta\times[-1,\infty)
\cap \mu_X^{-1}\bigl([k_\rho,\infty)\bigr)$.

First form $M\ecs U-Z(\rho)$ and see that this is a smooth manifold
with boundary $\mu_X^{-1}(k_\rho)$ union 
$\tau\bigl(\partial\Delta\times[k_\rho,\infty)\bigr)$.
Let $X(\rho)=\bigl(M\ecs U-Z(\rho)\bigr)\cup
\mu_X^{-1}(k_\rho)\times[k_\rho,\infty)$.

Next check that $M(\rho)\subset X(\rho)$ and that any dual to $w_2$
in $M(\rho)$ will be a dual to $w_2$ in $X(\rho)$.
This follows because topologically $X(\rho)$ can be 
obtained from $M(\rho)$ by adding a Spin manifold along an
embedded $\R^3$.
It now follows that 
$\gamma\bigl(M(\rho)\bigr)\leq \gamma\bigl(X(\rho)\bigr)$.

Next check that $X(\rho)$ has few essential $3${--}handles.\nl
For all $i\geq k_\rho$,
$\mu_X^{-1}\bigl([-4,i]\bigr)=W^\prime_i$ is a smooth submanifold
and $\ctrlBar{W^\prime_{i+1}-W^\prime_i}01$ is built topologically 
from $\ctrlBar{W_{i+1}-W_i}01$ by removing a copy of
$\Delta\times[i,i+1]$ and replacing it with a copy of
$\mu_X^{-1}(k_\rho)\times[i,i+1]$.
Since $H_3(W_{i+1},W_i;\Z)=0$, $H_3(W^\prime_{i+1},W^\prime_i;\Z)=0$.

Finally, use the Mayer{--}Vietoris sequence to see
$\beta_2\bigl(X(\rho)\bigr)$ is finite.

\newR{xx}
\proclaim \fxx \rm
Some sort of finiteness is needed to make 
the proof of \fmsofM\ work. 
If $M$ is the connected sum of infinitely many copies of 
$S^2\times S^2$, all the $M(s)$ constructed above are diffeomorphic.
To see this note that there is a smooth proper 
$h${--}cobordism, $H$, between $\R^4$ and $\open U(s)$.
Embed a ray cross $[0,1]$ properly into $H$ and remove a tubular
neighborhood.
Do the same with $M\times[0,1]$ and glue the two manifolds 
together along the boundary $\R^3\times[0,1]$.
The result is a smooth proper $h${--}cobordism between $M$ and $M(s)$.
This proper $h${--}cobordism is trivial by the usual argument involving 
summing a proper $h${--}cobordism with sufficiently many 
$S^2\times S^2\times[0,1]$.
Note that we are not saying that $M$ has only one smoothing but 
any different smoothing will have to have only finitely
many smooth $S^2\vee S^2$'s.

F~Ding and \v Z~Bi\v zaca and J~Etnyre have obtained results
on existence of smoothings which overlap the above results.
If $M$ is the interior of a topological manifold with a compact 
boundary component such that this component smoothly embeds 
in $k$ copies of $CP^2$, then F.~Ding \cite[\bDing] shows that 
there exists an $S>0$ so that all the $M(s)$ are distinct for $s\geq S$.
Bi\v zaca and Etnyre \cite[\bBE] show there are at least countably many 
distinct smoothings of $M$ among the $M(s)$ if $M$ is a topological 
manifold with a compact boundary component.
The Bi\v zaca and Etnyre proof is clearly a precursor of the proof above.

An example not covered by either
Ding or Bi\v zaca and Etnyre can be constructed as follows.
Use \fegAnyA\ to construct an $M$ with few $\ft${--}essential 
$3${--}handles such that $H_1(M;\Z)=\Q$ and $H_2(M;\Z)=0$.
Then $M$ has at least countably many distinct smoothings, possibly
the first such example with one end which is not topologically collared.

\newsec{\number\sectionnumber\quad Infinite Covers}

In the theorems below, we will use the hypothesis that $M$ is an infinite 
cover to bound 
$\displaystyle\lim_{r\to\infty}\gamma(\mecs r{\R_e})$ 
for $e\in\calE(M)$.
All the examples for which we know this limit satisfy
$\displaystyle\lim_{r\to\infty}\gamma(\mecs r{\R_e})=0$ or $\infty$.

\newT{fza}
\proclaim \ffza
Let $M$ be a non-compact smooth Spin $4${--}manifold 
with an action by $\Z$ which is smooth and properly discontinuous.
Then, for any $e\in\calE(M)$ and any integer $0\leq r<\infty$, 
$\gamma(\mecs r{\insm e})\leq\gamma(M)$.
If $\gamma(M)<\infty$ then for some
$e\in\calE(M)$, 
$\displaystyle\lim_{r\to\infty}\gamma(\mecs r{\insm e}) 
= \gamma(M)$.

\pf
Properly discontinuous implies that the orbit space $N$ is a smooth 
manifold and that the map $M\to N$ is a covering projection.
Given any $e\in\calE(M)$, we can find some 
element $g\in \Z$ such that all the translates of $e(D^4)$ 
by the powers of $g$ are disjoint, so all the 
$g^i\circ e\in\calE(M)$. 
For any integer $r>0$, use \fwd\ to find $e^\prime\in\calE(M)$ such that 
$\open D^4_{e^\prime}=\mecs r{\open D^4_e}$. 
The result follows.\break\hbox{} \qed

\newR{xx}
\proclaim \fxx \rm
Let $M$ be a smooth manifold homeomorphic to $\R^4$.
If $M$ is a non{-}trivial cover of some other $4${--}manifold 
then $M$ has a properly{-}discontin\-uous $\Z$ action.
Gompf \Cite[\bKirb, Problem 4.79A] asks for such smoothings and
remarks that most smoothings of $\R^4$ are not universal covers
of compact manifolds because there are only countably many 
smooth compact manifolds but there are uncountably many
smoothings of $\R^4$.
\abrev{goktso}\ gives the first concrete examples of smoothings
of $\R^4$ which can not be universal covers.

\newT{rffg}
\proclaim \frffg
Let $N$ be a compact smooth Spin $4${--}manifold with an infinite, 
residually-finite fundamental group.
Let $M$ denote the universal cover and let ${\cal C}$
denote the set of finite sheeted covers.
Then 
$$\gamma(\mecs r{\insm e})\leq
\gamma(M)=\max_{P\in{\cal C}}\{\gamma(P)\}$$
for all $e\in\calE(M)$.

\pf
Since there are elements of infinite order in $\pi_1(M)$,
\ffza\ implies $\gamma(\mecs r{\insm e})\leq
\gamma(M)$.
Let $e\in\calE(M)$.
It follows from the residual finiteness of $\pi_1(N)$ that 
there is a finite cover $P$ so that the composite
$D^4\ \RA{e}\ M\to P$ represents an element in
$\calE(P)$.
Hence $\gamma(M)\leq \max_{P\in{\cal C}}\{\gamma(P)\}$.
Given any $e\in\calE(P)$, $e$ lifts to an infinite number of
disjoint copies in $M$, so 
$\max_{P\in{\cal C}}\{\gamma(P)\}\leq\gamma(M)$.\break
\hbox{}\qed

\newE{xx}
\proclaim \fxx \rm
If $N$ is a smoothing of $S^3\times S^1$ 
then $\gamma(M)=0$.
If $N$ is a smoothing of the four torus
then $\gamma(M)\leq3$.
To see these results just observe that all the finite covers, $P$, 
of $N$ are homeomorphic to $N$.
These $P$ are compact, Spin and have signature $0$.
Hence $\gamma(P)=0$ if $N=S^3\times S^1$ 
and $\gamma(P)\leq 3$ for $N$ the four torus. 
Now apply \frffg.

\np
\goodbreak\vskip 15pt plus10pt minus5pt 
{\large\bf References}\ppar     
\leftskip=25pt\frenchspacing    
\parskip=3pt plus2pt\small      
\def\,{\thinspace}\let\sl\it

\bibitem{bizaca.z.&entyre.j.1996.1}%
{\bf{\v Z}~Bi{\v z}aca}, {\bf J~Etnyre}, {\it Smooth structures on collarable ends of
  4{--}manifolds}, preprint, 1996

\bibitem{casson.a.j.1986.1}%
{\bf A\,J~Casson}, {\it Three lectures on new{-}infinite constructions in
  $4${--}dimensional manifolds}, Progress in Math. (Birkhauser, Boston)
  62 (1986) 201--244 

\bibitem{demichelis.s&freedman.m.h.1992.1}%
{\bf S~DeMichelis},  {\bf M\,H~Freedman}, {\it Uncountably many exotic ${\R}^4$'s in
  standard $4${--}space}, J. Differential Geometry 35 (1992) 219--254

\bibitem{ding.f.1997.1}%
{\bf F~Ding}, {\it Uncountably many smooth structures on some 
  4{--}manifolds}, preprint (1996)

\bibitem{donaldson.s.k.1986.1}%
{\bf S\,K~Donaldson}, {\it Connections, cohomology and the intersection 
 forms of $4${--}manifolds}, J. Differential Geometry {24} (1986) 275--341

\bibitem{freedman.m.h.1979.1}%
{\bf M\,H Freedman}, {\it A fake {$S^3\times\R$}}, Ann. of Math. {110} 
 (1979) 177--201

\bibitem{freedman.m.h.1982.1}%
{\bf M\,H Freedman}, {\it The topology of four{--}dimensional manifolds}, 
  J. Differential Geometry {17} (1982) 357--453

\bibitem{freedman.m.h.&taylor.l.r.1986.1}%
{\bf M\,H~Freedman},  {\bf L\,R~Taylor}, {\it A universal smoothing of four{--}space},
  J. Differential Geometry {24} (1986) 69--78

\bibitem{furuta.m.1995.1}%
{\bf M~Furuta}, {\sl Monopole equation and the $11/8${--}conjecture}, 
  preprint RIMS, Kyoto (1995)

\bibitem{gompf.r.1993.1}%
{\bf R~Gompf}, {\sl An exotic menagerie}, J. Differential Geometry {37} 
  (1993) 199--22.

\bibitem{gompf.r.1996.1}%
{\bf R~Gompf}, {\sl Handlebody construction of {S}tein surfaces}, (1996)
  \hfill\penalty-1000\hbox{\tt ftp://ftp.ma.utexas.edu/pub/papers%
/gompf/stein\_FINAL.ps} 

\bibitem{gompf.r.1997.1}%
{\bf R~Gompf}, {\sl Kirby calculus for {S}tein surfaces}, (1997)
  \hfill\penalty-1000\hbox{\tt ftp://ftp.ma.utexas.edu/pub/papers/%
combs/stein-abridged.ps} 

\bibitem{grauert.h&remmert.r.1979.1}%
{\bf H~Grauert}, {\bf R~Remmert}, {\it {T}heory of {S}tein {S}paces}, 
  Grundl. Math. Wiss. vol. 226, Springer{--}Verlag, Berlin (1979)

\bibitem{kirby.r.c.1989.1}%
{\bf R\,C~Kirby}, {\it {T}opology of $4${--}{M}anifolds}, 
  Lecture Notes in Math. vol. 1374, Sprin\-ger{--}Verlag, Berlin (1989)

\bibitem{kirby.r.c.1996.1}%
{\bf R\,C~Kirby}, {\sl {P}roblems in {L}ow-{D}imensional {T}opology}, 
  ``{G}eometric {T}opology, {P}roc. of the 1993 {G}eorgia {I}nternational 
  {T}opology {C}onference'' (W\,H Kazez, ed), vol. 2 Part 2, Amer. Math. 
  Soc. and International Press (1996) 35--473

\bibitem{kirby.r.c&siebenmann.l.c.1977.1}%
{\bf R\,C~Kirby}, {\bf L\,C~Siebenmann}, {\it {F}oundational {E}ssays on 
  {T}opological {Manifolds}, {S}moothings, and {T}riangulations}, Annals of 
  Math. Studies, vol.~88, Princeton University Press, Princeton (1977)

\bibitem{milnor.j.1963.1}%
{\bf J~Milnor}, {\it Morse Theory}, Annals of Math. Studies, vol.~51, 
  Princeton University Press, Princeton (1963)

\bibitem{taubes.c.h.1987.1}%
{\bf C\,H~Taubes}, {\sl Gauge theory on asymptotically periodic 
  $4${--}manifolds}, J. Differential Geometry {25} (1987) 363--430

\bibitem{wall.c.t.c.1964.1}%
{\bf C\,T\,C~Wall}, {\it Diffeomorphisms of $4${--}manifolds}, 
  J. London Math. Soc. {39} (1964) 131--140

\end